\documentclass{article}
\usepackage[cp1251]{inputenc}
\usepackage[russian]{babel}
\usepackage{amsmath,amsthm}
\usepackage{amsfonts}
\usepackage{cite}
\def\g{\ensuremath{\mathfrak g}}
\def\G{\ensuremath{\mathbb G}}

\newtheorem{thm}{Теорема}
\newtheorem{lem}{Лемма}

\theoremstyle{definition}

\theoremstyle{remark}
\newtheorem{rem}[thm]{Замечание}
\DeclareMathOperator\adm{adm}
\DeclareMathOperator\Adm{Adm}

\title{Равенство емкости  и модуля конденсатора в субфинслеровом пространстве}
\author{Ю. В. Дымченко}
\date{}
\begin{document}
\maketitle
Равенство емкости и модуля конденсатора имеет важное значение в геометрической теории
функций. Оно позволяет связать теоретико-функциональные и геометрические свойства множеств. Для конформных емкостей и модулей на плоскости равенство было доказано Л. Альфорсом и А. Бёрлингом в работе \cite{ahlfors}. Затем этот результат был улучшен в работах Б. Фюгледе \cite{fuglede} и В. Цимера \cite{ziemer}. Дж. Хессе \cite{hesse} распространил этот результат на p-емкость и p-модуль
для случая, когда пластины конденсатора не пересекаются с границей области. В случае евклидовой метрики равенство емкости и модуля в самых общих предположениях было доказано В. А. Шлыком \cite{shlyk2}, затем это доказательство было немного упрощено в работе М. Оцука \cite{ohtsuka}. В случае римановой метрики равенство было доказано в \cite{capmod}.

Финслеровы пространства были введены как обобщение римановых многообразий на случай, когда метрика зависит не только от координат, но и от направления.
Равенство емкости и модуля конденсатора в финслеровых пространствах в самых общих предположениях было установлено в работе \cite{dymch2009}.

Пространства Карно-Каратеодори и субфинслеровы пространства отличаются от римановых и финслеровых пространств соответственно ограничением класса допустимых путей. C основными вопросами  анализа на группах Карно можно ознакомиться, например, в книге \cite{folland}. Емкости, модули конденсаторов, а также свойства различных функциональных классов на группах Карно в последнее время изучались  группой С. К. Водопьянова (например, \cite{vod89,vod98,vod96}). В частности, равенство емкости и модуля конденсатора  было установлено И. Г. Маркиной в работе \cite{markina2003}.

Субфинслеровы пространства изучались, например, в работах \cite{clelland2006628,ber13,donne,ber14,buk14}.

Приведем основные определения и обозначения. Доказательство многих нижеприведенных рассуждений можно найти в \cite{folland}.

Стратифицированной однородной группой (или группой Карно) называется связная односвязная нильпотентная группа Ли \G, алгебра Ли которой \g{} разлагается в прямую сумму векторных пространств $V_1\oplus V_2\oplus\dots\oplus V_m$ таких, что $[V_1,V_k]=V_{k+1}$ для $k=1,2,\dots,m-1$ и $[V_1,V_m]=\{0\}$. Здесь $[X,Y]=XY-YX$ --- коммутатор элементов $X$ и $Y$, а $[V_1,V_j]$ --- линейная оболочка элементов $[X,Y]$, где $X\in V_1$, $Y\in V_j$, $j=1,2,\dots,m$.

Пусть левоинвариантные векторные поля $X_{11}$, $X_{12}$,\dots, $X_{1n_1}$ образуют базис $V_1$. Определим подрасслоение $HT$ касательного расслоения $T\G$ со слоями $HT_x$, $x\in \G$, которые представляют собой линейную оболочку векторных полей $X_{11}(x)$, $X_{12}$, \dots, $X_{1n_1}(x)$. Назовем $HT$ горизонтальным касательным расслоением, а его слои $HT_x$ --- горизонтальными касательными пространствами в точке $x\in\G$.

Расширим базис $X_{11}$, \dots, $X_{1n_1}$ до базиса $X_{ij}$, $j=1,2,\dots n_i$, $i=1,2,\dots,m$, всей алгебры Ли \g, где каждый $X_{ij}$ представляет собой коммутатор $j$-го порядка некоторых векторов $X_{1j}$, $j=1,2,\dots,n_1$. Таким образом, $n_i$ является размерностью пространства $V_i$, $i=1,2,\dots,m$.

Любой элемент $x\in \G$ можно единственным образом представить в виде $x=\exp\left(\sum\limits_{i,j}x_{ij}X_{ij}\right)$. Набор чисел $\{x_{ij}\}$ назовем координатами элемента $x$. Получим взаимно однозначное отображение между группой \G{} и пространством $R^N$, где $N=n_1+n_2+\cdots+n_m$ --- топологическая размерность группы \G.

Мера Лебега в $R^N$ индуцирует биинвариантную меру Хаара в \G, которую мы обозначим через $dx$.

Обозначим $x_i=(x_{i1},x_{i2},\dots,x_{in_i})$, $i=1,2,\dots,m$. Определим растяжения $\delta_\lambda x$, $\lambda>0$, по формуле $\delta_\lambda x=(\lambda x_1,\lambda^2x_2,\dots,\lambda^m x_m)$. также имеем $d(\delta_\lambda x)=\lambda^Q dx$, где $Q=\sum\limits_iin_i$ --- однородная размерность группы \G.

Пусть  $F(x,\xi)$ --- неотрицательная функция, определенная при $x\in\G$, $\xi\in HT_x$,  которая гладко зависит от $x$ и $\xi$ и представляет собой финслерову метрику на каждом слое $HT_x$, то есть:

1) Для любого $a>0$ выполнено $F(x,a\xi)=aF(x,\xi)$ и $F(x,\xi)>0$ при $\xi\ne0$, $x\in\G$;

2) Для любых $x\in\G$, $\xi,\eta\in   HT_x$ функция $\nabla^2_HF^2(x,\eta)(\xi,\xi)$ положительно определена, где
$$
(\nabla^2_H)_{ij}=\frac12(X_{1i}X_{1j}+X_{1j}X_{1i}),\quad i,j=1,2,\dots,n_1.
$$

Определим на кокасательном расслоении $HT^*$ функцию $H(x,\omega)$, где $x\in\G$, $\omega\in HT^*_x$ как супремум величин $\omega(\xi)$ по всем $\xi\in HT_x$, удовлетворяющим условию $F(x,\xi)\le1$.  В дальнейшем будем отождествлять $\omega$ с вектором, имеющем координаты дифференциальной формы $\omega$ в базисе $\omega_i$, двойственным к базису $X_{1i}$, то есть $\omega_i(X_{1j})=\delta_{ij}$ для $i,j=1,2,\dots,n_1$.

Кривую $\gamma:(a,b)\to\G$ назовем горизонтальной, если для почти всех $t\in (a,b)$ $\dot \gamma(t)\in HT_{\gamma(t)}$. Длину такой кривой определим как интеграл $l(\gamma)=\int\limits_a^bF(\gamma(t),\dot\gamma(t))dt$. Если длина конечна, то кривую назовем спрямляемой.

На группе \G{} определим однородную норму $|\cdot|$, удовлетворяющую условиям: для любого $x\in\G$ $|x|\ge0$ и $|x|=0$ только при $x=0$; $|x^{-1}|=|x|$, $|\delta_\lambda x|=\lambda |x|$. Определим шар с центром в точке $x\in\G$ радиуса $r>0$ следующим образом: $B(x,r)=\{y\in \G: |x^{-1}y|<r\}$. Заметим, что он является левым сдвигом шара $B(0,r)$, который в свою очередь является образом единичного шара $B(0,1)$ при растяжении $\delta_r$. Известно \cite{folland}, что существует константа $C$ такая, что для любых $x,y\in\G$

\begin{equation}\label{ner}
\left||xy|-|x|\right|\le C|y|\text{ при } |y|\le \frac{|x|}2.
\end{equation}

Меру $dx$ нормируем так, чтобы $|B(0,1)|=\int\limits_{B(0,1)}dx=1$. очевидно, что $|B(0,r)|=r^Q$. Посредством непрерывной положительной в \G {} функции $g(x)$ определим элемент объема $d\sigma=g(x)dx$.

Расстояние $d_c(x,y)$ между двумя точками $x,y\in\G${} определим как инфимум длин кривых, соединяющих $x$ и $y$.

Пусть $D$ --- область в $\G$ и $E_0,E_1\subset \bar D$ ---  замкнутые непересекающиеся множества. Тройку множеств $(E_0,E_1,D)$ назовем конденсатором.

Будем говорить, что кривая $\gamma:(a,b)\to D$ соединяет множества $E_0$ и $E_1$, если $\liminf\limits_{t\to a} d(\gamma(t), E_0)=\liminf\limits_{t\to b} d(\gamma(t),E_1)=0$, где $d(x,y)=|x^{-1}y|$ для $x,y\in\G$. Семейство всех таких локально спрямляемых кривых обозначим через $\Gamma(E_0,E_1,D)$.

Расстояния $d$ и $d_c$ эквивалентны друг другу, а топология, порожденная расстоянием $d$, эквивалентна евклидовой \cite{vod98}.

Неотрицательную числовую борелевскую функцию на $D$ назовём допустимой для некоторого семейства $\Gamma$ кривых, расположенных в $D$,  если для любой $\gamma\in\Gamma$ $\int\limits_\gamma\rho F(x,dx)=\int\limits_a^b\rho(\gamma(t)) F(\gamma(t),\dot\gamma(t))dt \ge1$, где $\gamma(t)$ --- параметризация $\gamma$ посредством параметра $t\in(a,b)$. Множество всех допустимых функций для $\Gamma$ обозначим через $\adm\Gamma$.

Пусть $p>1$. Определим $p$-модуль конденсатора $(E_0,E_1,D)$ следующим образом:
$$
M_{p,F}(E_0,E_1,D)=\inf\int\limits_D\rho^p\,d\sigma,
$$
где инфимум берется по всем $\rho\in\adm\Gamma(E_0,E_1,D)$.

Функцию $u:D\to \mathbb R$ назовем локально липшицевой в $D$, если для любого компактного подмножества $D'\subset D$ существует константа $L$ такая, что для любых $x,y\in D'$ $u(x)-u(y)\le L d_c(x,y)$. Определим класс $L^1_{p,F}(D)$ как замыкание класса локально липшицевых в $D$ функций по норме
$$
\|u\|_{L^1_{p,F}(D)}=\left(\int\limits_DH(x,Xu)^p\,d\sigma\right)^{1/p},
$$
где $Xu=(X_{11}u,X_{12}u,\dots, X_{1n_1}u)$ --- горизонтальный градиент функции $u$, который имеет смысл в силу теоремы Радемахера для групп Карно (см. \cite{mitchell1985}).

Обозначим через $\Adm(E_0,E_1,D)$ множество неотрицательных функций из $L_{p,F}^1(D)\cap C(D)$, равных нулю (единице) в некоторой окрестности $E_0$ ($E_1$). Определим $p$-емкость конденсатора:
$$
C_{p,F}(E_0,E_1,D)=\inf\int\limits_DH(x,Xu)^p\,d\sigma,
$$
где инфимум берется по всем функциям $u\in\Adm(E_0,E_1,D)$.

\begin{lem}
  Инфимум в определении $M_{p,F}(E_0,E_1,D)$ можно брать по непрерывным в $D\setminus(E_0\cup E_1)$ допустимым функциям.
\end{lem}
\textbf{Доказательство.}

Пусть $0<\varepsilon<1/2$, $D_k$, $k=1,2,\ldots$  --- открытые множества, образующие
исчерпание изнутри множества $D\setminus( E_0\cup E_1)$, т.е. $\overline{D_k}\subset D_{k+1}$, $\bigcup\limits_{k=1}^\infty D_k=D\setminus( E_0\cup E_1)$;$d_k=d(\partial D_k,\partial D_{k+1})$, $k\ge1$. Положим для единообразия рассуждений $d_{-1}=d_0=\infty$, $D_0=\emptyset$.

Для каждого $k\ge1$ покроем компактное множество $\overline{D_k}\setminus D_{k-1}$ конечным числом шаров $B(x_i,r_i)$, где $x_i\in  \overline{D_k}\setminus D_{k-1}$, $r_i<\min(d_{k-2},d_k)/2$. Получим локально конечное покрытие области $D\setminus( E_0\cup E_1)$ шарами $B(x_i,r_i)$, $i\ge1$, лежащими в $D\setminus( E_0\cup E_1)$. Заметим, что покрытие шарами с теми же центрами и вдвое большими радиусами $B(x_i,2r_i)$ обладает тем же свойством. Дополнительно можно считать что все $r_i<1/2$.

Пусть $\{h_i(x)\}$ --- разбиение единицы на $D\setminus( E_0\cup E_1)$, подчиненное покрытию $\{B(x_i,r_i)\}$.

Возьмём допустимую функцию $\rho$ для $\Gamma(E_0,E_1,D)$ такую, что
$$
\int\limits_D\rho^p\,d\sigma<M_{p,F}(E_0,E_1,D)+\varepsilon.
$$
Пусть $\varphi(z)$ --- бесконечно дифференцируемая неотрицательная  в \G{} функция с носителем в $B(0,1)$ с условием $\int\limits_\G\varphi(z)\,dz=1$.

 Обозначим $\rho_i=h_i\rho$, $\varphi_t(x)=t^{-Q}\varphi(\delta_{1/t}x)$, $\tilde \rho_i=\int\limits_\G\rho_i(y)\varphi_t(xy^{-1})dy$. Для каждого $i\ge1$ подберем параметр $0<t_i<\varepsilon$ так, чтобы при $t\le t_i$ $\|\tilde\rho_i-\rho_i\|_{p,F}<2^{-i}\varepsilon^{1/p}$, где норма берется в пространстве $L_{p,F}(D)$, $\|\rho\|_{p,F}=\left(\int\limits_D\rho^p\,d\sigma\right)^{1/p}$ (см. \cite[утв. 1.20]{folland}, с заменой левых сдвигов на правые и наоборот). Также потребуем, чтобы $zB(x_i,r_i)\subset B(x_i,2r_i)$ для любого $z$ с $|z|\le t_i$. Это можно сделать в силу неравенства \eqref{ner}.

Функция $\log F(x,\xi)$ равномерно непрерывна  на компакте $\{(x,\xi):x\in \overline{B(x_i,2r_i)}, 1/2\le F(x,\xi)\le3/2\}$, то есть существует $\delta>0$ такое, что при $|z|<\delta$, $|\xi'-\xi''|<\delta$ и для любого $x\in B(x_i,2r_i)$ такого, что $zx\in \overline{B(x_i,2r_i)}$,
\begin{equation}
  \label{1}
  \frac{F(zx,\xi')}{F(x,\xi'')}\ge(1+\varepsilon)^{-1}.
\end{equation}

Здесь $\xi'$, $\xi''$ рассматриваем как векторы в базисе $X_{1i}$, $i=1,2,\dots,n_1$ с евклидовой нормой. Далее считаем, что $t_i<\delta$.

Функция $\tilde\rho=\sum\limits_i\tilde\rho_i$ является бесконечно дифференцируемой в $D\setminus( E_0\cup E_1)$ и

\begin{equation}\label{5}
\int\limits_D\tilde\rho^p\,d\sigma<M_{p,F}(E_0,E_1,D)+2\varepsilon,
\end{equation}

если положить $\tilde\rho=0$ на $E_0\cup E_1$.

Далее мы покажем, что функция $(1+\varepsilon)\tilde\rho$ является допустимой для $\Gamma(E_0,E_1,D)$.

Если $\gamma\in\Gamma(E_0,E_1,D)$, то
\begin{equation}
\label{3}
1\le\int\limits_\gamma\rho\,F(x,dx)=\int\limits_\gamma\sum_i\rho_i\,F(x,dx)=\sum_i\int\limits_{\gamma\cap B(x_i,r_i)}\rho_i\,F(x,dx).
\end{equation}

Преобразуем интеграл от функции $\tilde\rho$ по $\gamma$:
\begin{multline}
\label{2}
    \int\limits_\gamma\tilde\rho\,F(x,dx)=\int\limits_\gamma\sum_i\tilde\rho_i\,F(x,dx)=\int\limits_\gamma\sum_i\int\limits_\G\rho_i(y^{-1}x)\varphi_{t_i}(y)\,dy\,F(x,dx)= \\
   =\int\limits_\gamma\sum_i\int\limits_\G\rho_i((\delta_{t_i}z)^{-1}x)\varphi(z)\,dz\,F(x,dx)=
   \sum_i\int\limits_{B(0,1)}\varphi(z)\,dz\int\limits_{\gamma\cap B(x_i,2r_i)}\rho_i((\delta_{t_i}z)^{-1}x)\,F(x,dx)=\\=\int\limits_{B(0,1)}\varphi(z)\,dz\sum_i\int\limits_{\gamma\cap B(x_i,2r_i)}\rho_i((\delta_{t_i}z)^{-1}x)\,F(x,dx).
  \end{multline}

Пусть $z\in B(0,1)$. Для любого $i=1,2,\dots$ рассмотрим дугу $\gamma'$ из множества $\gamma\cap B(x_i,2r_i)$. Обозначим $\tilde\gamma'=(\delta_{t_i}z)^{-1}\cdot\gamma'$. Эта кривая будет горизонтальной вследствие левоинвариантности $X_{1j}$, $j=1,2,\dots,n_1$.  Соединим соответствующие граничные точки $\gamma'$ и $\tilde\gamma'$ двумя горизонтальными кривыми, лежащими в $B(x_i,2r_i)\setminus B(x_i,r_i)$. Заменим $\gamma'$ на объединение $\tilde\gamma'$ с этими горизонтальными кривыми. В результате всех этих изменений получим кривую $\tilde\gamma_z$, которая также будет допустимой для $\Gamma(E_0,E_1,D)$.

Далее имеем, делая замену $y=(\delta_{t_i}z)^{-1}x$ и параметризуя кривую $\gamma$ посредством финслеровой длины дуги (если кривая не спрямляема, то отсчитываем длину дуги от какой-либо точки кривой с соответствующим знаком):
$$
\int\limits_{\gamma\cap B(x_i,2r_i)}\rho_i((\delta_{t_i}z)^{-1}x)\,F(x,dx)\ge(1+\varepsilon)^{-1}\int\limits_{\tilde\gamma_z\cap B(x_i,r_i)}\rho_i(y)\,F(y,dy)
$$
в силу \eqref{1}. Подставляя в \eqref{2} и используя \eqref{3}, получим, что $\int\limits_\gamma\tilde\rho\,F(x,dx)\ge(1+\varepsilon)^{-1}$, то есть $(1+\varepsilon)\tilde\rho\in\adm \Gamma(E_0,E_1,D)$.

В силу произвольности $\varepsilon$ и неравенства \eqref{5} лемма доказана.

Следующий результат был установлен в $R^n$ В. А. Шлыком в \cite{shlyk2} и  модифицирован в работе \cite{ohtsuka1999}.

Определим систему замкнутых множеств $E_{ij}$, $j\ge0$, $i=0,1$, таких, что $E_{ij}\subset  \mathop{\rm int}E_{i,j-1}$ при $j\ge1$, $E_i=\bigcap\limits_{j=0}^\infty E_{ij}$, $E_{00}\cap E_{10}=\emptyset$.

\begin{lem}
\label{lemma}
  Пусть $\rho\in L_{p,F}(D)$ --- положительная непрерывная в $D\setminus(E_0\cup E_1)$ функция. Для любого $\varepsilon>0$ существует функция $\rho'$, $\rho'\ge\rho$ в $D$ такая, что
  \begin{enumerate}
    \item $\int\limits_D\rho'^p\,d\sigma\le\int\limits_D\rho^p\,d\sigma+\varepsilon$,
    \item Предположим, что для каждого $j\ge0$ существует кривая $\gamma_j\in\Gamma(E_{0j},E_{1j},D)$ такая, что\\ $\int\limits_{\gamma_j}\rho'\,F(x,dx)\le\alpha$. Тогда существует кривая $\tilde\gamma\in \Gamma(E_0,E_1,D)$ такая, что $\int\limits_{\tilde\gamma}\rho\, F(x,dx)\le\alpha+\varepsilon$.
  \end{enumerate}
  \end{lem}
\textbf{Доказательство. }

Пусть $E^j=E_{0j}\cup E_{1j}$, $W_j=E^{j-1}\setminus \mathop{\rm int}E^j$, $d_j=\min(d_c(\partial E_{0j},\partial E_{0,j-1}),d_c(\partial E_{1,j-1},\partial E_{1j}))>0$. Так как функция $\rho$ положительна в $D\setminus (E_0\cup E_1)$, можно найти последовательность $\varepsilon_j\to0$ при $j\to\infty$ такую, что
\begin{gather}
  \label{11}\sum_{j=1}^\infty(1+\varepsilon_j^{-1})\varepsilon_j^{p+1}<\varepsilon,\\
  \label{12}\alpha\varepsilon_j<d_j\inf_{W_j\cap D}\rho.
\end{gather}

Образуем последовательность компактных множеств $D_j$ такую, что $D_j\subset \mathop{\rm int} D_{j+1}$, $\bigcup\limits_{j=1}^\infty D_j=D$ и $\int\limits_{D\setminus D_j}\rho^p\,d\sigma<\varepsilon_j$.

Пусть $V_j=(D\setminus D_j)\cap W_j$. Положим
$$
\rho'(x)=\left\{
           \begin{array}{cl}
             (1+\varepsilon_j^{-1})\rho(x), & x\in V_j; \\
             \rho(x), & x\in D\setminus\bigcup\limits_{j=1}^\infty V_j.
           \end{array}
         \right.
$$
Покажем, что функция $\rho'$ удовлетворяет условиям леммы. Используя \eqref{11}, имеем:

\begin{align*}
  \int\limits_D\rho'^p\,d\sigma&=\sum_{j=1}^\infty\int\limits_{V_j}\left((1+\varepsilon_j^{-1})\rho\right)^p\,d\sigma+\int\limits_{D\setminus\bigcup\limits_{j=1}^\infty V_j}\rho^p\,d\sigma\le \\
  &\le \sum_{j=1}^\infty(1+\varepsilon_j^{-1})^p\int\limits_{V_j}\rho^p\,d\sigma+\int\limits_D\rho^p\,d\sigma\le \\
  &\le \sum_{j=1}^\infty(1+\varepsilon_j^{-1})^p\varepsilon_j^{p+1}+\int\limits_D\rho^p\,d\sigma\le\int\limits_D\rho^p\,d\sigma+\varepsilon.
\end{align*}
Таким образом, условие 1 выполнено. Покажем, что выполняется условие 2. Зафиксируем $j\ge1$. Кривая $\gamma_k\in\Gamma(E_{0j},E_{1j},D)$ при $k\ge j$, следовательно, она содержит две дуги: $\gamma'_k$, соединяющую $\partial E_{0j}$ с $\partial E_{0,j-1}$, и  $\gamma''_k$, соединяющую $\partial E_{1,j-1}$ с $\partial E_{1j}$. Дуги $\gamma'_k$ и $\gamma''_k$ не содержатся в $V_j$. Действительно, если бы выполнялось обратное, то с помощью неравенства \eqref{12} выводим, что
$$
\alpha\ge\int\limits_{\gamma_k}\rho'\,F(x,dx)\ge\int\limits_{\gamma'_k}\rho'\,F(x,dx)\ge\varepsilon_j^{-1}\int\limits_{\gamma'_k}\rho\,F(x,dx)\ge\varepsilon_j^{-1}d_j\inf_{W_j\cap D}\rho>\alpha,
$$
и аналогично с $\gamma''_k$. Получили противоречие. Значит,
$$
\gamma_k\cap(D_j\cap(E_{i,j-1}\setminus \mathop{\rm int}E_{ij}))\ne\emptyset,\quad i=0,1,\quad k\ge j.
$$
Обозначим $\gamma_k=\gamma_{0k}$. Приведем алгоритм, позволяющий из некоторой последовательности кривых $\gamma_{j-1,k}$ извлечь подпоследовательность $\gamma_{jk}$.

Заметим, что множество $D_j\cap(E_{i,j-1}\setminus \mathop{\rm int}E_{ij})$ является компактом. Следовательно, из последовательности $\gamma_{j-1,k}$ можно выделить подпоследовательность (которую снова обозначим $\gamma_{j-1,k}$), сходящуюся к  некоторой кривой $\gamma_0$, для которой множество $M=\gamma_0\cap(D_j\cap(E_{0,j-1}\setminus \mathop{\rm int}E_{0j}))\ne\emptyset$.
Возьмём какую-либо точку $x_{0j}\in M$. Так как $\rho$ непрерывна в точке $x_{0j}$, можно выбрать шар $B(x_{0j},r(x_{0j}))$ такой, что для любой геодезической линии $l$, соединяющую центр шара и его границу, выполнено условие
\begin{equation}
  \label{13}
\int\limits_l\rho\,F(x,dx)\le\frac\varepsilon{2^{j+3}}.
\end{equation}
Отбрасывая несколько первых членов  последовательности $\gamma_{j-1,k}$, можно считать, что любая кривая этой подпоследовательности пересекает шар $B(x_{0j},r(x_{0j}))$. Таким же образом рассмотрим множество $D_j\cap(E_{1,j-1}\setminus \mathop{\rm int}E_{1j})$, точку $x_{1j}$ из этого множества, шар $B(x_{1j},r(x_{1j}))$, удовлетворяющий условию, аналогичному \eqref{13}; так же из $\gamma_{j-1,k}$ выделим подпоследовательность, все члены которой пересекают этот шар. Полученная подпоследовательность и будет искомой последовательностью $\gamma_{jk}$.

Проводим изложенное построение последовательно для $j=1,2,\dots$. Рассмотрим диагональную последовательность $\gamma_{kk}$.

Кривая $\gamma_{kk}$ пересекает шары $B(x_{ij},r(x_{ij}))$, $i=0,1$, для $1\le j\le k$ не менее чем в двух точках. Соединим две точки пересечения с центром соответствующего шара геодезическими линиями. Получим кривую $\tilde\gamma_k\in\Gamma(E_{0k},E_{1k},D)$, проходящую через точки $x_{0j}$, $x_{1j}$, $j=1,2,\dots,k$. Для этой кривой имеем, используя условие \eqref{13}:
$$
\int\limits_{\tilde\gamma_k}\rho\,F(x,dx)\le\int\limits_{\gamma_{kk}}\rho\,F(x,dx)+2\sum_{j=1}^k\frac\varepsilon{2^{j+3}}\le\alpha+\frac\varepsilon4.
$$
Пусть $\Gamma_0$ --- семейство горизонтальных кривых, соединяющих $x_{00}$ и $x_{10}$ в $D\setminus(E_0\cup E_1)$, $\Gamma_{ij}$ --- семейство горизонтальных кривых в $D\setminus(E_0\cup E_1)$, соединяющих $x_{ij}$ и $x_{i,j+1}$, $i=0,1$, $j=1,2,\dots$.  Тогда
$$
\inf_{\gamma\in\Gamma_0}\int\limits_\gamma\rho\,F(x,dx)+\sum_{j=1}^k\inf_{\gamma\in\Gamma_{0j}}\int\limits_\gamma\rho\,F(x,dx)+\sum_{j=1}^k\inf_{\gamma\in\Gamma_{1j}}\int\limits_\gamma\rho\,F(x,dx)
\le\int\limits_{\tilde\gamma_k}\rho\,F(x,dx)\le\alpha+\frac\varepsilon4.
$$
Это верно для любого $k$, следовательно,
$$
\inf_{\gamma\in\Gamma_0}\int\limits_\gamma\rho\,F(x,dx)+\sum_{j=1}^\infty\inf_{\gamma\in\Gamma_{0j}}\int\limits_\gamma\rho\,F(x,dx)+\sum_{j=1}^\infty\inf_{\gamma\in\Gamma_{1j}}\int\limits_\gamma\rho\,F(x,dx)
\le\alpha+\frac\varepsilon4.
$$
Выберем кривые $C_0\in\Gamma_0$ и $C_{ij}\in\Gamma_{ij}$, $i=0,1$, $j=1,2,\dots$ так, чтобы

\begin{align*}
  \int\limits_{C_0}\rho\,F(x,dx)&<\inf_{\gamma\in\Gamma_0}\int\limits_\gamma\rho\,F(x,dx)+\frac\varepsilon2,\\
  \int\limits_{C_{ij}}\rho\,F(x,dx)&<\inf_{\gamma\in\Gamma_{ij}}\int\limits_\gamma\rho\,F(x,dx)+\frac\varepsilon{2^{j+3}}.
\end{align*}

Пусть $\tilde\gamma=\dots+C_{01}+C_0+C_{11}+\dots$. Тогда $\tilde\gamma\in\Gamma(E_0,E_1,D)$ и
$$
\int\limits_{\tilde\gamma}\rho\,F(x,dx)\le\alpha+\frac\varepsilon4+\frac\varepsilon2+2\sum_{j=1}^\infty\frac\varepsilon{2^{j+3}}=\alpha+\varepsilon.
$$
Лемма доказана.

Докажем теперь основную теорему.
\begin{thm}
  Пусть $D$ --- область в \G; $E_0$, $E_1$ --- непересекающиеся непустые компакты из $\bar D$. Тогда
  $$
  M_{p,F}(\Gamma(E_0,E_1,D))=C_{p,F}(E_0,E_1,D).
  $$
\end{thm}
\textbf{Доказательство. }
Сначала докажем неравенство
\begin{equation}\label{14}
  M_{p,F}(\Gamma(E_0,E_1,D))\le C_{p,F}(E_0,E_1,D).
\end{equation}
Пусть $u\in\Adm(E_0,E_1,D)$, $\Gamma_0$ --- подсемейство локально спрямляемых горизонтальных кривых $\gamma$ из $\Gamma(E_0,E_1,D)$ таких, что $u$ абсолютно непрерывна на любой спрямляемой замкнутой части $\gamma$. Определим функцию $\rho(x)=H(x,Xu)$ на $D$.

Пусть $\gamma\in\Gamma_0$ и $\gamma:(a,b)\to D$. Если $a<t_1<t_2<b$, то получим:
$$
\int\limits_\gamma\rho\,F(x,dx)\ge\int\limits_{t_1}^{t_2} H(x,Xu(\gamma(t)))F(x,\dot\gamma(t))dt\ge\left|\int\limits_{t_1}^{t_2} (Xu(\gamma(t)),\dot\gamma(t))\,dt\right|=|u(\gamma(t_2))-u(\gamma(t_1))|.
$$
Вследствие произвольности $t_1$ и $t_2$ получим, что $\int\limits_\gamma\rho\,F(x,dx)\ge1$. Таким образом, $\rho\in\adm\Gamma_0$.

Следовательно,
$$
M_{p,F}(\Gamma_0)\le\int\limits_D\rho^p\,d\sigma=\int\limits_DH(x,Xu)\,d\sigma.
$$
Учитывая, что $M_{p,F}(\Gamma_0)=M_{p,F}(\Gamma(E_0,E_1,D))$ (см. \cite{fuglede} и \cite{markina2003}), переходя к инфимуму по $u$, получим неравенство \eqref{14}.

Докажем противоположное неравенство
\begin{equation}\label{16}
  M_{p,F}(\Gamma(E_0,E_1,D))\ge C_{p,F}(E_0,E_1,D)
\end{equation}
в случае $(E_0\cup E_1)\cap\partial D=\emptyset$. Пусть $\rho\in \adm\Gamma(E_0,E_1,D)$ --- непрерывная в $D\setminus(E_0\cup E_1)$ функция. Определим в $D$ функцию $u(x)=\min(1,\inf\int\limits_{\beta_x}\rho\,F(x,dx))$, где инфимум берется по всем локально спрямляемым горизонтальным кривым $\beta_x$, соединяющим $E_0$ и $x$ в направлении точки $x$. Покажем, что $u\in\Adm(E_0,E_1,D)$ и $H(x,Xu)\le \rho$ почти везде в $D$. Если $u\equiv1$, то это очевидно.

В случае $u\not\equiv1$ пусть $\alpha_{x_1x_2}$ --- кратчайшая кривая, соединяющая $x_1$  и $x_2$ в направлении точки $x_2$, где точки $x_1$ и $x_2$ выбраны достаточно близко друг от друга. Пусть $\beta_{x_1}$ --- спрямляемая кривая, соединяющая $x_1$ и $E_0$. Тогда
$$
u(x_2)\le\int\limits_{\beta_{x_1}}\rho\,F(x,dx)+\int\limits_{\alpha_{x_1x_2}}\rho\,F(x,dx)\le\int\limits_{\beta_{x_1}}\rho\,F(x,dx)+\max_{x\in\alpha_{x_1x_2}}\rho(x)d_c(x_1,x_2).
$$
Здесь $d_c$ измеряется от $x_1$ до $x_2$.
Так как $\beta_{x_1}$ произвольно, то
$$
u(x_2)\le u(x_1)+\max_{x\in\alpha_{x_1x_2}}\rho(x)d_c(x_1,x_2).
$$
Используя рассуждения, аналогичные приведеным в \cite{mitchell1985}, убедимся в существовании производных $X_{1j}u$, $j=1,2,\dots,n_1$ почти всюду в $D$. Пусть $x_1$ такая точка, и пусть дана гладкая кривая, проходящая через $x_1$ в направлении вектора $\xi$.  Устремляя $x_2$ к $x_1$ по этой кривой, получим:
$$
Xu(x_1)(\xi)\le \rho(x_1)F(x_1,\xi).
$$
Поделив на $F(x_1,\xi)$ и взяв супремум по всем $\xi$, получим, что  $H(x_1,Xu(x_1))\le\rho(x_1)$. Значит,
$$
C_{p,F}(E_0,E_1,D)\le\int\limits_DH(x,Xu)^p\,d\sigma\le\int\limits_D\rho^p\,d\sigma.
$$
Переходя к инфимуму по $\rho$, получим неравенство \eqref{16} в случае $\partial D\cap(E_0\cup E_1)=\emptyset$. Поэтому теорема в этом случае доказана.

Рассмотрим общий случай $\partial D\cap(E_0\cup E_1)\ne\emptyset$. Пусть $0<\varepsilon<1/2$. Рассмотрим непрерывную на $D\setminus(E_0\cup E_1)$ допустимую для $\Gamma(E_0,E_1,D)$ функцию $\rho$ такую, что
$$
\int\limits_{D\setminus(E_0\cup E_1)}\rho^p\,d\sigma<\varepsilon+M_{p,F}(E_0,E_1,D).
$$
Можем считать, что $\rho>0$ на $D\setminus(E_0\cup E_1)$, иначе возьмем вместо нее функцию $\max(\rho(x),h(x))$, где $h(x)>0$ --- непрерывная на \G{} функция  со сколь угодно малым интегралом $\int\limits_Dh^p\,d\sigma$.

Пусть $\rho'$, $E_{0j}$, $E_{1j}$ такие же, как в лемме \ref{lemma}. Покажем, что
$$
\int\limits_\gamma\rho'\,F(x,dx)>1-2\varepsilon
$$
для всех $\gamma\in\Gamma(E_{0j},E_{1j},D)$ при достаточно больших $j$.

Действительно, если это не так, то найдутся $j_k$ и $\gamma_k\in\Gamma(E_{0j_k},E_{1j_k},D)$ такие, что $\int\limits_{\gamma_k}\rho'\,F(x,dx)\le1-2\varepsilon$. По лемме \ref{lemma} найдется кривая $\tilde\gamma\in\Gamma(E_0,E_1,D)$ такая, что $\int\limits_{\tilde\gamma}\rho\,F(x,dx)\le1-\varepsilon$, что противоречит допустимости функции $\rho$.

Определим функцию
$$
\tilde\rho(x)=\left\{
  \begin{array}{ll}
    \dfrac{\rho'}{1-2\varepsilon}, & x\in D\setminus(E_{0j}\cup E_{1j}); \\
    0, & x\notin D\setminus(E_{0j}\cup E_{1j}).
  \end{array}
\right.
$$
Она принадлежит $\adm\Gamma(E_0,E_1,D\cup E_{0j}\cup E_{1j})$. Поэтому, в силу доказанного частного случая,
\begin{multline*}
 C_{p,F}(E_0,E_1,D)\le C_{p,F}(E_0,E_1,D\cup E_{0j}\cup E_{1j})=M_{p,F}(E_0,E_1,D\cup E_{0j}\cup E_{1j})\le\\
\le\int\limits_D\tilde\rho^p\,d\sigma\le(M_{p,F}(E_0,E_1,D)+2\varepsilon)(1-2\varepsilon)^{-p}.
\end{multline*}

Устремляя $\varepsilon\to0$, получим неравенство \eqref{16}, а следовательно, и утверждение теоремы.
\begin{rem}
  Из леммы \ref{lemma} следует непрерывность модуля, то есть
$$
\lim\limits_{j\to\infty}M_{p,F}(E_{0j},E_{1j},D)=M_{p,F}(E_0,E_1,D).
$$
\end{rem}

%\bibliographystyle{pdmi}
%\bibliography{bibl}
\end{document}